\documentclass[10pt]{amsart}

\usepackage{amsfonts}
\usepackage{amsmath}

\begin{document}

\newcommand{\e}{\mathrm{E}}
\newcommand{\var}{\mathrm{Var}}
\newcommand{\deq}{\stackrel{d}{=}}
\newcommand{\dx}{\mathrm{d}}
\newcommand{\rdp}{\rho^{\prime\prime}}
\newcommand{\rdpx}{\rho^{\prime\prime}(x)}

\title{Symmetrization of Bernoulli}
\author{Soumik Pal \\ Department of Statistics, Columbia University}
\date{\today}
\maketitle

\section{Introduction}\label{intro}

Let $X$ be a random variable. We shall call an independent random variable $Y$ to be a {\emph {symmetrizer}} for $X$, if $X+Y$ is symmetric around zero. If $Y$ is independent copy of $-X$, it is obviously a symmetrizer. A random variable is said to be \emph{symmetry resistant} if the variance of any symmetrizer $Y$, is never smaller than the variance of $X$ itself. 

For example, let $X$ be a Bernoulli($p$) random variable. If $p=1/2$, it is immediate that the degenerate random variable, $Y\equiv -1/2$, is a symmetrizer for $X$. Hence, $X$ is not symmetry resistant. However, we shall show that if $p\neq 1/2$, for any symmetrizer $Y$, we have
\[
\var(Y) \geq pq,
\]
where $q = (1-p)$. It is immediate from this inequality that $X$ is symmetry resistant and the minimum variance symmetrizer has the same variance as $-X$.

The last result is the main content of a paper by Kagan, Mallows, Shepp, Vanderbei and Vardi (see \cite{symm}), where the reader can look for the motivation and the history of this problem. We merely reprove the result here. However, the novelty is that while the previous authors have used technology to solve what they see as {\lq}infinite dimensional linear programming problem{\rq}, we have used purely probabilitic techniques, as ubiquitous as Ito's rule. In fact, at a crucial point in the paper~\cite{symm}, the authors had to use a linear programming software; to quote: {\lq}We then used the fourth authorÕs software, called LOQO, to solve these discrete approximations and eventually were able to guess at the correct functional form for $\rho$.{\rq} This was to the discomfort of Dr. Larry Shepp (private communications), one of the authors, to whom I am thankful for suggesting this problem.

As said before, our proof is completely probabilistic, essentially relying on Skorokhod embedding and Ito's rule for Stochatic Calculus.

\section{Proof}\label{proof}
Let $Y$ be any symmetrizer for $X$ with finite variance. By Skorokhod embedding, there is a stopping time $\tau$ such that for any standard Brownian motion (i.e., starting from zero), the stopped process has the distribution of $Y - \e(Y)$. That is to say, if $W$ is a standard Brownian motion
\begin{equation}\label{skem}
 W_{\tau} \deq Y - \e(Y).
\end{equation}
Here and throughout $\deq$ refers to equality in distribution. 
\smallskip

On a suitable probability space construct a process $\{B_t\}$ such that
\begin{equation}\label{b0}
B_0 \deq X - \e(X)
\end{equation}
and
\[
B_t = B_0 + W_t,
\]
where $B_0$ is independent of the standard Brownian motion $W_t$. Then clearly, $B_t$ is a Brownian motion which has the initial distribution of $X - \e(X)$. Also, by equation~(\ref{skem}), we have
\[
B_{\tau} \deq X + Y - \e(X+Y).
\]
But, since $Y$ is a symmetrizer of $X$, we should have $\e(X+Y)=0$, and hence
\begin{equation}\label{skem2}
B_{\tau} \deq X + Y.
\end{equation}
\smallskip

Now, let $\rho$ be any smooth odd function with bounded derivatives on the real line. By Ito's rule, we have
\begin{equation}\label{ito}
\rho(B_{t}) - \rho(B_0) = M_t + \frac{1}{2}\int_0^{t}\rdp(B_s)\dx s,
\end{equation}
where $M_t$ is a martingale. Thus, by the optional sampling theorem, we have
\begin{equation}\label{ito2}
\e(\rho(B_{\tau})) - \e(\rho(B_0)) = \frac{1}{2} \e\left(\int_0^{\tau}\rdp(B_s)\dx s\right).
\end{equation}

Now, since $\rho$ is an odd function, and $X+Y$ is symmetric around zero, by equation~(\ref{skem2}), $\e(\rho(B_{\tau})) = 0$. Thus~(\ref{ito2}) reduces to
\begin{equation}\label{ito3}
- \e(\rho(B_0)) = \frac{1}{2} \e\left(\int_0^{\tau}\rdp(B_s)\dx s\right).
\end{equation}
Let us now look at the RHS of ~(\ref{ito3}). By conditioning on $B_0$, we have
\begin{eqnarray}\label{rhs}
\e\left(\int_0^{\tau}\rdp(B_s)\dx s\right)&=& p\e\left(\int_0^{\tau}\rdp(q + W_s)\dx s\right)\nonumber\\
 &+& q\e\left(\int_0^{\tau}\rdp(-p + W_s)\dx s\right).
\end{eqnarray}
Let now impose the following restrictions on $\rho$ (an example will soon follow):

\begin{enumerate}
\item $\vert \rdp \vert \leq 1$.
\item $\rho(1 + x) = -\rho(x), \;\; \forall x \in {\bf R}$. 
\end{enumerate}
in addition to the fact that $\rho$ is odd. Thus for any $x$, we have \[
\rdp(-p + x) = -\rdp( 1-p + x) = -\rdp(q + x).
\]
Then, by equation ~(\ref{rhs}) we get
\[
\e\left(\int_0^{\tau}\rdp(B_s)\dx s\right)=(p-q)\e\left(\int_0^{\tau}\rdp(q + W_s)\dx s\right).
\]
Also, by equation ~(\ref{b0}),
\[
\e(\rho(B_0)) = p\rho(q) + q\rho(-p) = p\rho(q) - q\rho(-p) = (p-q)\rho(q).
\] 
Substituting these values in equation ~(\ref{ito2}) will give
\begin{equation}\label{main}
(q-p)\rho(q) = -\frac{1}{2}(q-p)\e\left(\int_0^{\tau}\rdp(q + W_s)\dx s\right).
\end{equation}
\smallskip

\noindent Without loss of generality, assume that $q < p$. The case when $q = p$ can be easily dealt separately. The other case of $q > p$ simply needs some sign flippings in the following argument. Now, use the fact that $\vert \rdp \vert \leq 1$, to conclude
\[
-\e\left(\int_0^{\tau}\rdp(q + W_s)\dx s\right) \leq \e(\tau).
\]
Thus, from ~(\ref{main}), we get
\[
\e(\tau) \geq 2\rho(q).
\]
However, by Skorokhod embedding, $\e(\tau) = \var(Y)$. Hence
\begin{equation}\label{main2}
\var(Y) \geq 2\rho(q).
\end{equation}
\smallskip

\noindent Finally, we have to exhibit such a $\rho$. This is how we define it. For $0 \leq x \leq 1$, define
\[
\rho(x) = \frac{x(1-x)}{2}
\]
Extend it to the entire positive axis by the property $\rho(1+x)=-\rho(x)$. That is to say,
\[
\rho(x) := (-1)^{k}\rho( x - k ), \;\; 0 \leq k \leq x \leq k+1 .
\]
And extend to the entire negative axis by the oddness of $\rho$. That is
\[
\rho(-x) := - \rho(x), \;\; x \in {\bf R}^{+}.
\] 
The function $\rho$ does not have a continuous second derivative. However, the set of discontinuity is just the countable set of integers, and this is sufficient for the usual Ito's rule to go through. See, for example,~\cite{bmsc}, page 219. More importantly, 
\[
\vert \rdpx \vert \leq 1,
\]
whenever $x$ is not an integer. Hence by equation ~(\ref{main2}), we get
\[
\var(Y) \geq pq,
\]
which is what we claimed. This proves the result.\hfill$\Box$
\medskip
\hrule
\bigskip

\noindent It is interesting to see how every inequality above becomes equality only when $\tau$ is the Skorokhod embedding for $Y = - X$. This proves the uniqueness of the minimum variance symmetrizer.


\begin{thebibliography}{80}
\bibitem{symm} {\sc{Kagan, A., Mallows, C., Shepp, L., Vanderbei, R. \& Vardi, Y.}} (1999) Symmetrization of Binary Random Variables. \emph{Bernoulli}, {\bf 5}, 1013-1020.
\bibitem{bmsc} {\sc{Karatzas, I. \& Shreve, S. E.}} (1991) \emph{Brownian Motion and Stochastic Calculus}. Second Edition, Springer-Verlag GTM, New York.



\end{thebibliography}
\end{document}